\documentclass{article}%
\usepackage{amsmath}
\usepackage{amssymb}%
\usepackage{amsfonts}%
\usepackage{graphicx}
\providecommand{\U}[1]{\protect\rule{.1in}{.1in}}
\newtheorem{theorem}{Theorem}

\newtheorem{definition}[theorem]{Definition}

\begin{document}

\title{A Gliding Hump Characterization of Normed Barrelled Spaces}
\author{Christopher Stuart\\Department of Mathematical Sciences\\New Mexico State University\\Las Cruces, New Mexico, 88003\\cstuart@nmsu.edu}
\maketitle

\begin{abstract}
In this paper, a characterization of normed barrelled spaces is given, as well
as a similar characterization of rings of sets with the Nikodym Property.
Finally, a condition that is equivalent to a Banach space having a separable
quotient is discussed.

\end{abstract}

\textbf{Introduction}

\qquad The Uniform Boundedness Principle (UBP) is one of the fundamental
principles in functional analysis. The locally convex topological vector
spaces for which the equicontinuity version of the UBP holds are known as
barrelled spaces. Recall that a barrelled space is one in which every barrel
(an absolutely convex, closed, absorbing set) is a neighborhood of zero. It is
not hard to show that in a normed linear space $X$, equicontinuity is
equivalent to uniform boundedness of a sequence of pointwise bounded
functionals $\left(  y_{n}\right)  \subset X^{\ast}.$ See [10] for the basic
results on these spaces.

\qquad To show a particular space is barrelled, a method of proof known as a
gliding or sliding hump argument is frequently used. In this paper, we present
a gliding hump characterization of barrelledness on a normed linear space $X$,
and also of the Nikodym Property on rings of sets. We also discuss the
relationship between non-barrelledness and the separable quotient problem in
Banach spaces.

The "gliding hump method" is an inelegant term for an elegant technique.
Basically, given a norm unbounded $\left(  y_{n}\right)  \subset X^{\ast},$ we
find a bounded sequence $\left(  x_{n}\right)  \subset X$ for which
$\left\vert y_{n}\left(  x_{n}\right)  \right\vert \rightarrow\infty$ as
$n\rightarrow\infty.$ We can think of $\left(  x_{n}\right)  $ as a
"detecting" sequence since it detects the unboundedness of $\left(
y_{n}\right)  .$ We then "collect" a subsequence $\left(  x_{n_{k}}\right)  $
of $\left(  x_{n}\right)  $ in an $x\in X$ for which $\left\vert y_{n_{k}%
}\left(  x\right)  \right\vert \rightarrow\infty$ as $k\rightarrow\infty.$ The
"hump" in the method is the sequence $\left\vert y_{n_{k}}\left(  x_{n_{k}%
}\right)  \right\vert $ which is chosen to have most of the "mass" of $\left(
y_{n_{k}}\right)  $ on $x.$ From this we can conclude that every norm
unbounded $\left(  y_{n}\right)  $ is pointwise unbounded, which establishes
the barrelledness of $X,$ since we will have proved the contrapositive. See
[1], [7], [8], [9] for many applications of this method.\bigskip

\textbf{Main Results}

\begin{definition}
Let $X$ be a normed linear space. A detecting sequence in $X$ for an unbounded
$\left(  y_{n}\right)  \subset X^{\ast}$ is a bounded sequence $\left(
x_{n}\right)  $ such that $\left\vert y_{n}\left(  x_{n}\right)  \right\vert
\rightarrow\infty$ as $n\rightarrow\infty.$
\end{definition}

In previous publications, the term "bounding sequence" was used instead of
"detecting sequence." The referee of this paper pointed out that, since the
sequence $\left(  x_{n}\right)  $ is actually "unbounding", a different term
seemed appropriate.

\begin{theorem}
\bigskip A normed linear space $X$ is barrelled if and only if for every norm
unbounded $\left(  y_{n}\right)  \subset X^{\ast}$ there is a detecting
sequence $\left(  x_{n}\right)  $ and an $x\in X$ such that for any
subsequence $\left(  x_{n_{k}}\right)  $ there is a further subsequence, also
denoted as $\left(  x_{n_{k}}\right)  ,$ for which we have $\left\vert
y_{n_{k}}\left(  x_{n_{k}}\right)  \right\vert -\left\vert y_{n_{k}}\left(
x-x_{n_{k}}\right)  \right\vert \rightarrow\infty$ as $k\rightarrow\infty. $
\end{theorem}

\textbf{Proof: }For sufficiency, we need to show that $\left(  y_{n_{k}%
}\left(  x\right)  \right)  $ is unbounded. But that is easy since
\[
\left\vert y_{n_{k}}\left(  x\right)  \right\vert \geq\left\vert y_{n_{k}%
}\left(  x_{n_{k}}\right)  \right\vert -\left\vert y_{n_{k}}\left(
x-x_{n_{k}}\right)  \right\vert \rightarrow\infty\text{ as }k\rightarrow\infty
\]
by the triangle inequality and the hypothesis in the theorem.

For necessity, assume that $X$ is barrelled. So for every norm unbounded
$\left(  y_{n}\right)  \subset X^{\ast}$ there is $x\in X$ with $\left\vert
y_{n}\left(  x\right)  \right\vert \rightarrow\infty$ as $k\rightarrow\infty.$
Choosing $x_{n}=x$ for all $n$ does satisfy the hypotheses in the theorem for
a suitable detecting sequence, but we are primarily interested in finding a
\emph{non-trivial }\ detecting sequence. That is, a sequence $\left(
x_{n}\right)  $ for which $\left\vert y_{n}\left(  x_{n}\right)  \right\vert
\rightarrow\infty$ but $\left(  y_{n}\right)  $ is pointwise bounded on
$\left\{  x_{n}\right\}  .$ This leads to the standard collection process for
a gliding hump argument, such as in the following, which is Theorem 2 in [9].

\begin{theorem}
$X$ is barrelled if for every subsequence $\left(  x_{n_{k}}\right)  $ of a
detecting sequence $\left(  x_{n}\right)  $ for an unbounded sequence $\left(
y_{n}\right)  \subset X^{\ast},$ there is a further subsequence $\left(
x_{n_{k_{l}}}\right)  $ such that $\sum_{l=1}^{\infty}x_{n_{k_{l}}}=x,$ where
the convergence is with respect to the weak topology on $X$ induced by
$\left(  y_{n}\right)  ,$ denoted as the $w\left(  y_{n}\right)  $ topology.
\end{theorem}

Within the proof of this theorem, the inequality%

\[
\left\vert y_{n_{k_{l}}}\left(  x\right)  \right\vert =\left\vert y_{n_{k_{l}%
}}\left(  \sum_{j=1}^{l-1}x_{n_{k_{j}}}\right)  +y_{n_{k_{l}}}\left(
x_{n_{k_{l}}}\right)  +y_{n_{k_{l}}}\left(  \sum_{j=l+1}^{\infty}x_{n_{k_{j}}%
}\right)  \right\vert
\]%
\[
\geq\left\vert y_{n_{k_{l}}}\left(  x_{n_{k_{l}}}\right)  \right\vert
-\left\vert y_{n_{k_{l}}}\left(  \sum_{j=1}^{l-1}x_{n_{k_{j}}}\right)
\right\vert -\left\vert y_{n_{k_{l}}}\left(  \sum_{j=l+1}^{\infty}x_{n_{k_{j}%
}}\right)  \right\vert \rightarrow\infty
\]
is derived, which is essentially the condition we need to satisfy the theorem.
$\square$

So it is the existence of a non-trivial detecting sequence that is the crux of
the gliding hump method. The Josefson-Nissenzweig Theorem guarantees that the
dual of any infinite-dimensional Banach space $X$ contains a sequence $\left(
y_{n}\right)  $ such that $\left\Vert y_{n}\right\Vert =1$ and $y_{n}\left(
x\right)  \rightarrow0$ as $n\rightarrow\infty$ for all $x\in X.$ We can find
a subsequence $\left(  y_{n_{k}}\right)  $ of $\left(  y_{n}\right)  $ and
$\left(  x_{n_{k}}\right)  \subset X,$ $\left\Vert x_{n_{k}}\right\Vert =1$
for all $k,$ such that $\left\vert y_{n_{k}}\left(  x_{n_{k}}\right)
\right\vert =1$ and $k^{2}\left\vert y_{n_{k}}\left(  x_{n_{l}}\right)
\right\vert \rightarrow0$ as $k\rightarrow\infty$ for all $l.$ So $\left(
ky_{n_{k}}\right)  $ is a norm unbounded sequence that is pointwise bounded on
$span\left\{  x_{n_{k}}\right\}  .$ There is no guarantee, though, that
$\left(  ky_{n_{k}}\right)  $ is pointwise bounded on a dense subspace of $X$,
which is needed to show that $X$ has an infinite-dimensional separable
quotient. In the last section of this paper, we will briefly discuss the
relationship of such sequences to the separable quotient problem in Banach spaces.

We next show that a condition that characterizes the Nikodym Property (NP) for
rings of subsets of a given set can be shown in a way similar to Theorem 2.
The original Nikodym Boundedness Theorem was proved for countably additive set
functions defined on $\sigma$-algebras of sets. It has been generalized to
finitely additive functions defined on rings of sets that satisfy various
completeness conditions. For a fairly general condition that implies (NP) see [7].

Recall that the elements of $ba\left(  R\right)  $ of are the bounded,
finitely additive measures on a ring $R.$

\begin{definition}
Let $\left(  y_{n}\right)  \subset ba\left(  R\right)  $ be an unbounded
sequence. A \emph{detecting sequence }for $\left(  y_{n}\right)  $ is a
sequence $\left(  E_{n}\right)  \subset R$ such that $\left\vert y_{n}\left(
E_{n}\right)  \right\vert \rightarrow\infty$ as $n\rightarrow\infty.$
\end{definition}

Note that $\left(  E_{n}\right)  $ is not assumed to be pairwise disjoint.

\begin{theorem}
A ring $R$ has NP if and only if for every unbounded $\left(  y_{n}\right)
\subset ba\left(  R\right)  $ there is a set $E\in R$ and a detecting sequence
of subsets of $E,$ $\left(  E_{n}\right)  ,$ such that for any subsequence
there is a further subsequence $\left(  E_{n_{k}}\right)  $ of $\left(
E_{n}\right)  $ we have
\[
\left\vert y_{n_{k}}\left(  E_{n_{k}}\right)  \right\vert -\left\vert
y_{n_{k}}\left(  E\backslash E_{n_{k}}\right)  \right\vert \rightarrow
\infty\text{ as }k\rightarrow\infty.
\]

\end{theorem}

\textbf{Proof: }To show sufficiency, we need to show that for $\left(
y_{n}\right)  $ unbounded on $R,$ we have a set $E\in R$ with $\left\vert
y_{n_{k}}\left(  E\right)  \right\vert \rightarrow\infty$ for some subsequence
of $\left(  y_{n}\right)  .$ This follows easily from the assumptions in the
theorem since%
\[
y_{n_{k}}\left(  E\right)  =y_{n_{k}}\left(  E_{n_{k}}\cup\left(  E\backslash
E_{n_{k}}\right)  \right)  =y_{n_{k}}\left(  E_{n_{k}}\right)  +y_{n_{k}%
}\left(  E\backslash E_{n_{k}}\right)
\]

so by the triangle inequality we have%

\[
\left\vert y_{n_{k}}\left(  E_{n_{k}}\right)  \right\vert =\left\vert
y_{n_{k}}(E)-y_{n_{k}}\left(  E\backslash E_{n_{k}}\right)  \right\vert
\leq\left\vert y_{n_{k}}\left(  E\right)  \right\vert +\left\vert y_{n_{k}%
}\left(  E\backslash E_{n_{k}}\right)  \right\vert .
\]
This gives%

\[
\left\vert y_{n_{k}}\left(  E\right)  \right\vert \geq\left\vert y_{n_{k}%
}\left(  E_{n_{k}}\right)  \right\vert -\left\vert y_{n_{k}}\left(
E\backslash E_{n_{k}}\right)  \right\vert \rightarrow\infty\text{ as
}k\rightarrow\infty.
\]

To show necessity, we assume that $R$ has the NP. So that for every unbounded
$\left(  y_{n}\right)  $ there is a set $E\in R$ such that $\left\vert
y_{n_{k}}\left(  E\right)  \right\vert \rightarrow\infty$ for some subsequence
of $\left(  y_{n}\right)  .$ It's true that setting $E_{n}=E$ would satisfy
the condition in the theorem, but we are primarily interested in finding a
pairwise disjoint detecting sequence for $\left(  y_{n}\right)  $ that
satisfies the condition. This is done in [9] with the following condition:

\begin{definition}
The ring $R$ is said to have the bounded subsequential completeness property
(BSCP) if for every unbounded $\left(  y_{n}\right)  \subset ba\left(
R\right)  $ there is a detecting sequence $\left(  E_{n}\right)  \subset R$
for which every subsequence $\left(  E_{n_{k}}\right)  $ has a further
subsequence $\left(  E_{n_{k_{l}}}\right)  $ with $\cup_{l}E_{n_{k_{l}}}\in
R.$
\end{definition}

The proof that a ring with the BSCP has NP is very similar to that given for a
normed barrelled space. The inequality derived satisfies the condition in this
theorem. $\square$

\bigskip

\textbf{Detecting sequences and dense, non-barrelled subspaces}

\bigskip

In this section, we are interested in the problem of constructing a dense
subspace $S$ of a Banach space $X$ that contains a non-trivial detecting
sequence for some norm unbounded $\left(  y_{n}\right)  \subset X^{\ast}.$ In
particular, we want $\left(  y_{n}\right)  $ to be pointwise bounded on $S.$
This is equivalent to $S$ being a non-barrelled, dense subspace of $X.$ Saxon
and Wilansky in [5] showed that $X$ having a non-barrelled, dense subspace is
equivalent to $X$ having a separable (infinite-dimensional) quotient.

It is not hard to see that the idea of a \emph{pseudo-orthogonal family }(see
S.A. Morris and D.T. Yost, [4]) gives a non-trivial detecting sequence. In
that paper, Morris and Yost show that a Banach space with a pseudo-orthogonal
family has a separable quotient. Also, Sliwa [6] introduces the idea of a
\emph{pseudobounded }sequence $\left(  y_{n}\right)  \subset X^{\ast},$ which
is a norm unbounded sequence \ $\left(  y_{n}\right)  $ that is pointwise
bounded on a dense subspace $S\subset X.$ As noted above, $S$ would then be a
dense, non-barrelled subspace of $X.$In Sliwa's paper, it is shown that $X$
having a dense, non-barrelled subspace is equivalent to $X^{\ast}$ containing
a $weak^{\ast}$ basic sequence $\left(  y_{n}\right)  .$ This, by Proposition
II.1 in [3], is equivalent to $X/\left\{  y_{n}\right\}  ^{\top}$ having a
basis $\left(  x_{n}\right)  .$ The notation $\left\{  y_{n}\right\}  ^{\top
}=\left\{  x\in X:y_{n}\left(  x\right)  =0\text{ for all }n\right\}  $ is
from [3].

Also as a result of Saxon and Wilansky's paper, we know that $X$ has a dense,
non-barrelled subspace if and only if $X$ has a proper quasicomplemented
subspace. This means that there exist infinite-dimensional subspaces $A$ and
$B$ of $X$ such that $A$ and $B$ are closed, $A\cap B=\left\{  0\right\}  ,$
and $A+B$ is dense in $X.$

This seems closely related to the existence of \emph{hereditarily
indecomposable }(H.I.) Banach spaces. A Banach space $X$ is H.I. if there does
not exist infinite-dimensional subspaces $A$ and $B$ of $X$ such that $X=A+B$
and $A\cap B=\left\{  0\right\}  ,$ and the same condition applies to any
infinite-dimensional close subspace of $X.$ See [2] for an introduction to
this topic.

\bigskip\newpage

\textbf{References}

[1] P. Antosik and C. Swartz, \emph{Matrix Methods in Analysis, }Lecture Notes
in Mathematics (Book 1113), Springer, 1985.

[2] H. H. Giv, \emph{What is a Hereditarily Indecomposable Banach Space?,
}Notices of the American Mathematical Society, Volume 66, No. 9, October 2019, 1474-1475.

[3] W.B. Johnson and H.P. Rosenthal, \emph{On }$\omega^{\ast}$\emph{-basic
sequences and their applications to the study of Banach spaces, }Studia
Mathematica, T. XLIII, 1972, 77-92.

[4] S.A. Morris and D.T. Yost, \emph{Observations of the Separable Quotient
Problem for Banach Spaces, }Axioms 2020 9, 7.

[5] S.A. Saxon and A. Wilansky \emph{The equivalence of some Banach space
problems, }Colloq. Math. 1977 37, 217-226.

[6] W. Sliwa, \emph{The separable quotient problem and the strongly normal
sequences, }J. Math. Soc. Japan, Vol. 64, No. 2 2012, 387-397.

[7] C.E. Stuart and P. Abraham, \emph{Generalizations of the Nikodym
Boundedness and Vitali-Hahn-Saks Theorems, }J. of Mathematical Analysis and
Applications 300 (2004), 351-361.

[8] C.E. Stuart, \emph{Normed Barrelled Spaces, }Algebraic Analysis and
Related Topics, Banach Center Publications, Volume 53, Institute of
Mathematics, Polish Academy of Sciences, Warszawa 2000, 205-210.

[9] C.E. Stuart, \emph{Normed Barrelled Spaces, }Proyecciones Journal of
Mathematics, Vol. 39, No. 5, October 2020, 1267-1272.

[10] C. Swartz, \emph{Introduction to Functional Analysis, }Marcel Dekker, 1992.

\end{document}